\documentclass{gtmon_a}
\pdfoutput=1

\proceedingstitle{Heegaard splittings of 3--manifolds (Haifa 2005)}
\conferencestart{10 July 2005}
\conferenceend{19 July 2005}
\conferencename{Heegaard splittings of 3--manifolds}
\conferencelocation{Haifa}

\editor{Cameron Gordon}
\givenname{Cameron}
\surname{Gordon}

\editor{Yoav}
\givenname{Yoav}
\surname{Moriah}

\title{Hyperbolic volume, Heegaard genus and ranks of groups}

\author{Peter B Shalen}
\givenname{Peter B}
\surname{Shalen}
\address{Department of Mathematics, Statistics, and Computer Science\\
University of Illinois at Chicago\\\newline
851 S. Morgan Street\\
Chicago IL 60607-7045\\
USA}
\email{shalen@math.uic.edu}
\urladdr{http://www.math.uic.edu/~shalen/}

\volumenumber{12}
\issuenumber{}
\publicationyear{2007}
\papernumber{14}
\startpage{335}
\endpage{349}

\doi{}
\MR{}
\Zbl{}

\arxivreference{}  

\keyword{hyperbolic manifold}
\keyword{volume}
\keyword{Heegaard genus}
\keyword{rank of group}
\subject{primary}{msc2000}{57M05}
\subject{primary}{msc2000}{57M50}

\received{29 September 2006}
\revised{27 January 2007}
\accepted{4 March 2007}
\published{3 December 2007}
\publishedonline{3 December 2007}
\proposed{}
\seconded{}
\corresponding{}
\version{}

\makeatletter
\def\cnewtheorem#1[#2]#3{\newtheorem{#1}{#3}[section]
\expandafter\let\csname c@#1\endcsname\c@para}

\let\xysavmatrix\xymatrix
\def\xymatrix{\disablesubscriptcorrection\xysavmatrix}
\AtBeginDocument{\let\bar\wbar\let\tilde\wtilde\let\hat\what}

\theoremstyle{definition}
\newtheorem{para}{\hspace{-2pt}}[section]
\cnewtheorem{remark}[para]{Remark}
\cnewtheorem{remarks}[para]{Remarks}
\cnewtheorem{definition}[para]{Definition}
\cnewtheorem{definitions}[para]{Definitions}
\cnewtheorem{notation}[para]{Notation}
\cnewtheorem{convention}[para]{Convention}

\theoremstyle{plain}
\cnewtheorem{theorem}[para]{Theorem}
\cnewtheorem{conjecture}[para]{Conjecture}
\cnewtheorem{lemma}[para]{Lemma}
\cnewtheorem{question}[para]{Question}
\cnewtheorem{proposition}[para]{Proposition}
\cnewtheorem{corollary}[para]{Corollary}
\cnewtheorem{idea}[para]{Idea}
\cnewtheorem{consequence}[para]{Consequence}

\numberwithin{equation}{para}

\def\dnewtheorem#1[#2]#3{\newtheorem{#1}{#3}[section]
\expandafter\let\csname c@#1\endcsname\c@thm}

\dnewtheorem{lem}[thm]{Lemma}

\theoremstyle{definition}
\dnewtheorem{defn}[thm]{Definition}

\makeatother  

\newcommand\Number{\begin{para}}
\newcommand\EndNumber{\end{para}}
\newcommand\Conjecture{\begin{conjecture}}
\newcommand\EndConjecture{\end{conjecture}}
\newcommand\Question{\begin{question}}
\newcommand\EndQuestion{\end{question}}
\newcommand\Theorem{\begin{theorem}}
\newcommand\EndTheorem{\end{theorem}}
\newcommand\Lemma{\begin{lemma}}
\newcommand\EndLemma{\end{lemma}}
\newcommand\Proposition{\begin{proposition}}
\newcommand\EndProposition{\end{proposition}}
\newcommand\Corollary{\begin{corollary}}
\newcommand\EndCorollary{\end{corollary}}
\newcommand\Proof{\begin{proof}}
\newcommand\EndProof{\end{proof}}
\newcommand\NoProof{{\hfill$\square$}}

\newcommand\WHAT{3.66}
\newcommand\barchi{\wbar\chi}
\newcommand\chibar{\barchi}
\newcommand\calw{{\mathcal W}}
\newcommand\Mthin{M_{\rm thin}}
\newcommand\Mthick{M_{\rm thick}}
\newcommand\Isom{{\rm Isom}}
\newcommand\ZZ{{\mathbb Z}}
\newcommand\HH{{\mathbb H}}
\newcommand\calh{{\mathcal H}}
\newcommand\vol{\mathop{\rm vol}}
\newcommand\tuberad{\mathop{{\rm tuberad}}}
\newcommand\drillc{\mathop{{\rm drill}_C}}
\newcommand\thebound{1.22}
\newcommand\rank{\mathop{{\rm rank}}}
\newcommand{\tM}{\wwtilde M}

\begin{document}

\begin{asciiabstract}
Some conjectures about Heegaard genera and ranks of fundamental groups
of 3-manifolds are formulated, and it is shown that they imply new
statements about hyperbolic volume.
\end{asciiabstract}

\begin{abstract} 
Some conjectures about Heegaard genera and ranks of fundamental groups
of 3--manifolds are formulated, and it is shown that they imply new
statements about hyperbolic volume.
\end{abstract}

\maketitle

\section{Preface}
In my talk at the conference on Heegaard splittings at the Technion, I
formulated some topological conjectures and
discussed how, by adapting some of my work with Marc Culler and others
on volumes of hyperbolic $3$--manifolds, to use these conjectures---if
true---to the problem of relating hyperbolic volume to Heegaard genus.

One of the conjectures that I stated is a 
modernized version of the antique conjecture that if $M$ is a compact,
orientable $3$--manifold, the rank of $\pi_1(M)$ is equal to the
Heegaard genus of $M$.  The first counterexamples to this old
conjecture, in which $M$ is a Seifert fibered manifold, had been given
by Boileau and Zieschang in \cite{bz}; more general counterexamples,
for graph manifolds, were given by Schultens and Weidmann in
\cite{schultweid}. Here is the modernized version:

\Conjecture\label{rank is genus}If $M$ is a compact, orientable,
hyperbolic $3$--manifold, the rank of $\pi_1(M)$ is equal to the
Heegaard genus of $M$.
\EndConjecture

The main results that I discussed in my
talk give connections between
Heegaard genus and hyperbolic volume that are conditional on
\fullref{rank is genus} and another conjecture, which is
formulated below as \fullref{genus version}. Those results are
incorporated into this article as Corollaries \ref{A'}, \ref{Ccor} and
\ref{BCor}. They are immediate conequences of results, \fullref{A} and Propositions \ref{C} and \ref{B}, which involve ranks of
groups instead of Heegaard genus and are not conditional on \fullref{rank is genus} (although they are conditional on other
topological conjectures). 

This reorganization of the material may prove valuable if \fullref{rank is genus} turns out to be false. The immediate motivation
for reorganizing the paper in this way was that rumors were
circulating to the effect that \fullref{rank is genus} had been
disproved. These rumors seem to have died down.

In this article I will provide detailed proofs of all the new results,
including those that I announce in my talk.  I will give less space here
than in my talk to summaries of the arguments in Agol--Culler--Shalen
\cite{rank4} and Culler--Shalen \cite{first}, on
which the proofs are based, as these summaries were not significantly
different from the ones given in the introductions to those papers.

All the new results are conditional in the sense that they include
various topological conjectures as hypotheses. Some of these
conjectures are very new, and I hope that they will motivate new
research, whether or not they are true.

I thank Ian Agol for some valuable discussions about \fullref{rank is genus}. I also thank the referee for pointing out some
major errors in the first draft of the article.

\section{Introduction}

When one studies closed hyperbolic $3$--manifolds, the volume is a key
invariant, because it is known that up to isometry there exist at most
finitely many hyperbolic manifolds of a given finite volume. The
volume of a closed hyperbolic $3$--manifold $M$ is a topological
invariant of $M$---for example because Mostow rigidity says the
hyperbolic structure of $M$ is itself determined by the topology---but
the precise connection between volume and more classical topological
invariants is far from being well understood.

In this article I'll be concerned with the connection between volume
of a hyperbolic $3$--manifold and the rank of its fundamental group.
An upper bound on rank---or even the Heegaard genus---of a closed
hyperbolic $3$--manifold does not give an upper bound on the
volume. For example, if $M_0$ is a closed hyperbolic $3$--manifold that fibres
over the circle, and $g$ is the genus of the fiber, then for every
positive integer $n$ there is an $n$--fold cyclic cover $M_n$ of $M$
which itself fibers over $S^1$ with genus--$g$ fiber. In particular the
Heegaard genus of each $M_n$ is at most $2g+1$, but the volume of
$M_n$ is $nv_0$, where $v_0$ denotes the volume of $M_0$.

On the other hand, an upper bound on the volume of a closed hyperbolic
$3$--manifold definitely does give an upper bound on the Heegaard genus
(and hence on the rank). In fact, there is a universal constant
$\mu>0$, the {\it Margulis constant,} such that for every closed
hyperbolic $3$--manifold $M$, the subset $\Mthin$ of $M$, defined to
consist of all points through which there pass homotopically
non-trivial curves of length $\le\mu$, is a disjoint union of smooth
solid tori. It's not hard to show that the submanifold
$\Mthick=\overline{M{-}\Mthin}$ has a triangulation with at most $C_0v$
simplices, where $v$ denotes the volume of $M$ and $C_0$ is another
universal constant that can be computed from $\mu$. This implies that
there is a constant $C$ such that every closed $3$--manifold $M$ has
Heegaard genus at most $Cv$, where $v$ denotes the volume of $M$.

The problem is that the constant $C$ that comes from this argument is
astronomical, and the estimates obtained in this way don't get us
anywhere near what we expect from examples. The goal of this article
will be to suggest a way of getting good explicit upper bounds for
the rank of $\pi_1(M)$, where $M$ is a hyperbolic $3$--manifold, in
terms of the volume of $M$.

By contrast, if one is willing to settle for bounds on the homological
complexity of a manifold $M$ instead of the rank of $\pi_1(M)$, the
results that appear in my papers with Culler and others are in the realm
of reality and are sometimes sharp. For example, \cite[Theorem~1.1]{rank4}
states among other things that if $M$ is a closed, orientable hyperbolic
$3$--manifold with volume at most $\thebound$, then $H_1(M;\ZZ_p)$ has
dimension at most $2$ for every prime $p\ne2,7$. This result is sharp
for $p=3$ and for $p=5$: the manifolds referred to in the Weeks--Hodgson
census \cite{census} as {\tt m003(-3,1)}         and {\tt m007(3,1)}
have respective volumes $0.94\ldots$    and $1.01\ldots$, while their
integer homology groups are respectively isomorphic to $\ZZ_5\oplus\ZZ_5$
and $\ZZ_3\oplus\ZZ_6$.

One would like to obtain bounds for Heegaard genus, or rank of
fundamental group, in terms of volume, similar to the bounds on
homology that I have described. This will require more knowledge on
the topological side. In \fullref{conjecture section} I'll
formulate and discuss a pair of topological conjectures, one about
rank and one about Heegaard genus; the former will be shown to imply a
bound on the rank of the fundamental group in terms of hyperbolic
volume. In \fullref{C-section} I'll show how to get a stronger
bound by combining the conjecture of \fullref{conjecture section}
with a conjectured analogue for rank of Moriah and Rubinstein's result
in \cite{mr} about the behavior of Heegaard genus under Dehn filling.
In \fullref{B-section} I'll give another consequence of the
conjecture of \fullref{conjecture section}, partially analogous to
the results of \cite{first} and \cite{second}.

\section{A little background}\label{baby section}

The following result provides the simplest motivation (from the point
of view of studying volumes) for the conjectures on ranks of
$3$--manifold groups that I'll state in the next section.

\Theorem\label{1.015}  Suppose that $M$ is a closed, orientable hyperbolic
$3$--manifold with volume less than $1.015$. Then $\pi_1(M)$ has a
$2$--generator subgroup of finite index.
\EndTheorem

\Proof Let $M$ be any closed, orientable hyperbolic $3$--manifold. Let
us write $M=\HH^3/\Gamma$, where $\Gamma\le\Isom^+(\HH^3)$ is discrete,
cocompact and torsion-free. If $\Gamma\cong\pi_1(M)$ has no $2$--generator
subgroup of finite index, it follows from Jaco and Shalen \cite[Theorem
VI.4.1]{JS} that $\Gamma$ is $2$--free in the sense that all its
$2$--generator subgroups are free.

When $\Gamma$ is $2$--free, the number $\log3$ is a strong Margulis number
for $M$ in the sense of Anderson, Canary, Culler and Shalen \cite[Section
8]{accs}. Indeed, the case $k=2$ of \cite[Proposition 8.1]{accs}
asserts that $\log3$ is a strong Margulis number for $M$ as long as
every $2$--generator subgroup of $\Gamma$ is free and topologically
tame. But according to the main results of Agol \cite{agoltameness} and
Calegari--Gabai \cite{cg}, every finitely generated Kleinian group is
topologically tame.

Przeworski \cite[Theorem~3]{prezbetti} asserts that if the first Betti number
$\beta_1(M)$ is at least $3$, then the volume of $M$ is at least
$1.015$. The hypothesis $\beta_1(M)\ge3$ is used only in order to
quote \cite[Theorem~1]{prezbetti}, which is included in 
\cite[Corollary~10.6]{accs}, and gives a lower bound $V(\lambda)$ on the
volume of $M$ in terms of the length $\lambda$ of a shortest geodesic
in $M$, under the assumption that $\beta_1(M)\ge3$. Now, according to
\cite[Corollary 10.5]{accs}, the number $V(\lambda)$ is a lower bound
for the volume of $M$ as long as $\log3$ is a strong Margulis number
for $M$.  Hence the lower bound of $1.015$ for the volume of $M$ still
holds if one replaces the assumption $\beta_1(M)\ge3$ by the
assumption that $\log3$ is a strong Margulis number for $M$; by the
discussion above, this holds in particular if $\pi_1(M)$ has no
finite-index subgroup of rank $2$. \EndProof

\section{A conjecture on rank and finite covering spaces, and a consequence}\label{conjecture section}

\fullref{1.015} raises the following question (most immediately
for the case $k=2$):

\Question  What restriction does the existence of  a
 $k$--generator subgroup of finite index in $\pi_1(M)$ place on $M$?
\EndQuestion

There is a {\it homological} condition 
which is necessary for the existence of such a subgroup: $H_1(M;{\bf
Z}_p)$ must have rank $\le k+1$ for every prime $p$.  The necessity of
this condition follows from Shalen--Wagreich \cite[Proposition
1.1]{shalen-wagreich}. This was the starting point for the results
relating volume to homology which I described in the introduction.

The following conjecture, which I first formulated in my talk at the
Technion, would provide a bound of the  rank of $\pi_1(M)$ when
$\pi_1(M)$ has a $k$--generator
subgroup of finite index.

\Conjecture\label{my conjecture}  If $M$ is a compact, orientable hyperbolic
$3$--manifold such that $\rank\pi_1(M)=r$, then for any finite-sheeted
covering space $\tM$ of $M$  we have $\rank\pi_1(\tM)\ge r-1$. 
\EndConjecture

Equivalently, \fullref{my conjecture} asserts that if 
$\pi_1(M)$ has a $k$--generator
subgroup of finite index then $\rank\pi_1(M)\le k+1$.

In my talk, I also formulated a parallel  conjecture about Heegaard
genus:  

\Conjecture\label{genus version}If $M$ is a compact, orientable hyperbolic
$3$--manifold with Heegaard genus $g$, then any finite-sheeted
covering space $\tM$ of $M$  has Heegaard genus at least $g-1$. 
\EndConjecture

Of course, Conjectures \ref{my conjecture} and \ref{genus version} are
equivalent modulo the ``rank equals genus'' \fullref{rank is
  genus}.

\Number\label{third remark} {\rm It follows from Agol, Culler and
Shalen \cite[Corollary 7.3]{rank4} (which is in turn a refinement of
Shalen and Wagreich \cite[Proposition 1.1]{shalen-wagreich}, a result
which I mentioned above) that if $M$ is a closed, orientable hyperbolic
$3$--manifold such that $H_1(M;\ZZ_p)$ has rank $r$ for a given prime
$p$, then for any finite-sheeted covering space $\wwtilde M$ of $M$,
the rank of $H_1(\wwtilde M;\ZZ_p)$ is at least $r-1$. We may regard
Conjectures \ref{my conjecture} and \ref{genus version} as analogues,
for the rank of the fundamental group and the Heegaard genus, of this
result about the rank of the mod--$p$ homology.}  \EndNumber

\Number\label{second remark} {\rm There appears to be a huge class of
examples, for every $g\ge3$, in which a compact, orientable hyperbolic
$3$--manifold $M$ with Heegaard genus $g$ has a finite-sheeted covering
space $\tM$ of Heegaard genus exactly $g-1$. Alan Reid gave the first
such example, for $g=3$, in his paper \cite{reid}, which was inspired
by an earlier, weaker version of \fullref{1.015} which appeared in
Culler and Shalen \cite{paradoxical}. 

During the conference, Hyam Rubinstein pointed out a systematic way of
constructing such examples.  Suppose that $M$ is a closed, orientable
hyperbolic $3$--manifold containing a closed, non-orientable surface
$F$, and that the complement in $M$ of the interior of a regular
neighborhood $N$ of $F$ is a handlebody $J$ of genus $g-1$.  Since $N$
is a twisted $I$--bundle over $F$, it's easy to construct a $2$--sheeted
covering space $q:\tM\to M$ to which $J$ lifts, and in which
$q^{-1}(N)$ is a trivial $I$--bundle. It follows that $\tM$ has a
Heegaard splitting in which the handlebodies are isotopic to the two
lifts of $J$; in particular this splitting has genus $g-1$. On the
other hand, if $A$ is a vertical arc in the twisted $I$--bundle $N$,
and $T$ is a regular neighborhood of $A$ relative to $N$, then $J\cup
T$ and $\overline{N{-}T}$ are genus--$g$ handlebodies, which define a
genus--$g$ Heegaard splitting of $M$. It appears that for $g\ge3$ the
``generic'' situation is that both the genus--$g$ Heegaard splitting of
$M$ and the genus--$(g{-}1)$ Heegaard splitting of $\tM$ are
minimal. 

These examples help show why the lower bound $g-1$ is natural in
\fullref{genus version}. I would guess that in many of these
examples one can also show that $\rank\pi_1(M)=g$ and $\rank\pi_1(M)=g-1$,
which would help show why the lower bound $r-1$ is natural in
\fullref{my conjecture}.}

\EndNumber

\Number\label{first remark} {\rm During the conference, Andrew Casson pointed
out that \fullref{my conjecture} is
trivial for a cyclic regular covering.  This is because if any group
$G$ has a normal subgroup $N$ with $G/N$ cyclic, the rank of $G$ can
obviously exceed the rank of $N$ by at most $1$.}  \EndNumber

\Number {\rm The best known result in the direction of \fullref{genus
version} seems to be the result of Rieck and Rubinstein \cite{rr}, which
gives a lower bound for the Heegaard genus of a two-sheeted covering of
$M$ in terms of the Heegaard genus of $M$.} \EndNumber

\fullref{1.015} now has the following immediate consequence:

\Corollary\label{A} If \fullref{my
  conjecture} is true, then for every closed, orientable hyperbolic
$3$--manifold with volume at most $1.015$ we have
$\rank\pi_1(M)\le3$.\NoProof \EndCorollary

Since Conjectures \ref{my conjecture} and \ref{genus version} are
equivalent modulo the ``rank equals genus'' \fullref{rank
    is genus}, we also get:

\Corollary \label{A'} If Conjectures \ref{rank is genus} and
 \ref{genus version} are true, then for every closed, orientable
 hyperbolic $3$--manifold with volume at most $1.015$ we have
 $\rank\pi_1(M)\le3$.\NoProof \EndCorollary

\section{A conjecture on rank and Dehn filling, and a consequence}\label{C-section}

It is a consequence of the result proved by Moriah and Rubinstein in
\cite{mr} that if $N$ is a hyperbolic $3$--manifold of finite volume
with exactly one cusp, and if $g$ denotes the Heegaard genus of the
compact core $\wwhat N$ of $N$, then infinitely many Dehn fillings
of $\wwhat N$ yield closed manifolds of Heegaard genus exactly $g$.
(This was re-proved by a purely topological argument in Rieck and
Sedgwick \cite{rs}.)

In view of Moriah and Rubinstein's result, the following conjecture
would follow immediately from the ``rank equals genus'' \fullref{rank is genus}.

\Conjecture\label{missing} Suppose that $N$ is a hyperbolic
$3$--manifold of finite volume with exactly one cusp, and set
$r=\rank\pi_1(N)$. Let $\wwhat N$ denote the compact core of $N$.
Then there is an infinite sequence $(M_i)$ of manifolds obtained by
distinct Dehn fillings of $\wwhat N$ such that each $\rank\pi_1(M_i)=r$
for each $i$.
\EndConjecture

The main result of \cite{rank4}, Theorem 1.1, states that if $M$ is a
closed, orientable hyperbolic $3$--manifold with volume at most
$\thebound$, then $H_1(M;\ZZ_p)$ has dimension at most $2$ for every
prime $p\ne2,7$, and that $H_1(M;\ZZ_2)$ and $H_1(M;\ZZ_7)$ have
dimension at most $3$. The result had originally been proved in a
weaker form, which states that if $M$ has volume at most $\thebound$,
then $H_1(M;\ZZ_p)$ has dimension at most $3$ for every prime $p$. In
this section I'll show how to adapt the proof of this weaker theorem
to prove an analogous result about the rank of the fundamental group,
modulo Conjectures  \ref{my conjecture} and \ref{missing}. This
is \fullref{C} below. I have not thought about whether one can
prove a Heegaard-genus analogue of the strong form of \cite[Theorem
1.1]{rank4} modulo Conjectures  \ref{my conjecture} and \ref{missing}. 

\Proposition\label{C} If Conjectures \ref{my conjecture} and
\ref{missing} are true, then for every closed, orientable hyperbolic
$3$--manifold $M$ with volume at most $1.22$, we have
$\rank\pi_1(M)\le3$.  \EndProposition

You'll notice that the information given by \fullref{C} modulo
Conjectures \ref{my conjecture} and \ref{missing} is strictly stronger
than the information given by \fullref{A} modulo \fullref{my conjecture}
alone. The proof of \fullref{C} uses a lot more mathematics than that
of \fullref{A}. In particular, \cite[Lemma~3.2]{rank4}, which is quoted
in the proof below, depends on a result from Agol, Dunfield, Storm and
Thurston \cite{ast} which in turn relies on Perelman's estimates for
the Ricci flow with surgeries.

\Proof[Proof of \fullref{C}]
As in \cite{rank4}, we shall say that a hyperbolic manifold $M$ is
{\it exceptional}  if every shortest geodesic in $M$ has tube radius at
most $(\log 3)/2$.

We first prove the proposition in the case where $M$ is
non-exceptional. In this case, by definition, there is a shortest
geodesic $C$ in $M$ with $R=\tuberad(C) > (\log3)/2$. We set
$N=\drillc (M)$. Let $\calh$ denote the maximal cusp neighborhood in
$N$.  Since $R > (\log3)/2$, \cite[Lemma 3.2]{rank4} implies that
$\vol \calh < \pi$.

Now assume that $\rank\pi_1(M)\ge4$. Set $r=\rank\pi_1(\wwhat
N)$. It is obvious that $r\ge \rank\pi_1(M)$, so in particular $r\ge4$.
\fullref{missing} implies that there is an infinite sequence $(M_i)$ of
manifolds obtained by distinct Dehn fillings of $\wwhat N$ such that
$\rank\pi_1(M_i)=r$ for each $i$.  Since $r\ge4$, \fullref{my conjecture}
implies that for each $i$, every finite-index subgroup of $\pi_1(M_i)$
has rank at least $3$. Since $\pi_1(M_i)$ has no $2$--generator subgroup
of finite index, it follows from Jaco and Shalen \cite[Theorem VI.4.1]{JS}
that $\pi_1(M_i)$ is $2$--free in the sense that all its $2$--generator
subgroups are free.  \cite[Lemma~4.3]{rank4} then implies that
$\vol \calh\ge\pi$, a contradiction.  This completes the proof in the
non-exceptional case.

We now turn to the case where $M$ is exceptional. If $M$ is isometric
to the manifold $\vol3$ discussed by Jones and Reid \cite{jones-reid},
then $M$ can be obtained by Dehn filling from a once-punctured
torus bundle over $S^1$, and hence $\rank\pi_1(M)\le3$. If $M$ is not
isometric to $\vol3$, then according to \cite[Proposition 7.1]{rank4},
the group $\pi_1(M)$ has a finite-index subgroup of rank at most $2$.
But if $\rank\pi_1(M)\ge4$, \fullref{my conjecture} implies that
every finite-index subgroup of $\pi_1(M)$ has rank at least $3$. Hence
$\rank\pi_1(M)\le3$ in this case as well.
\EndProof

I pointed out in \fullref{conjecture section} that Conjectures
\ref{my conjecture} and \ref{genus version} are equivalent modulo the
``rank equals genus'' \fullref{rank is genus}. I pointed out at
the beginning of the present section that \fullref{missing}
would follow immediately from \fullref{rank is genus}. Hence the following corollary follows from
\fullref{C}.

\Corollary\label{Ccor} If Conjectures \ref{rank is genus} and
\ref{genus version} are true, then for every closed, orientable hyperbolic
$3$--manifold $M$ with volume at most $1.22$, the Heegaard genus of $M$
is at most $3$.  \EndCorollary

\section {A hybrid consequence of the conjecture on covering spaces}\label{B-section}

In \cite{second}, Culler and I prove:

\Theorem  If $M$ is a closed, orientable hyperbolic
$3$--manifold with volume at most $3.08$, then $H_1(M;{\bf Z}_2)$ has
rank at most $7$. 
\EndTheorem

The weaker version that the rank is
at most $10$ is somewhat easier to prove, and is established in \cite{first}.

In this section I'll show how to adapt the proof of this weaker result
to get the following result. The information given by this result
modulo \fullref{my conjecture} is not exactly an analogue of
the result of \cite{first} involving the rank of the fundamental
group, but rather a hybrid result involving both the rank of the
fundamental group and the rank of the mod--$2$ first homology.

\Proposition\label{B} If \fullref{my
  conjecture} is true, then for every closed, orientable hyperbolic
$3$--manifold $M$ with volume at most $3.08$, either $\rank\pi_1(M)\le10$, or
$\dim H_1(M;\ZZ_2)\le4$.  \EndProposition

The proof of \fullref{B} will depend on refining a number of
the results proved in \cite{first}. For the rest of this section I
shall use the notation and definitions of \cite{first}, including the
definition of a book of $I$--bundles. As in \cite{first}, I shall
write $\chibar(X)=-\chi(X)$, where $X$ is any space homeomorphic to a
finite polyhedron and $\chi(X)$ denotes its Euler characteristic. 

\Lemma\label{new lemma}Let $X$ be a compact, connected $3$--manifold,
and let $P\ne X$ be a submanifold of $X$. Suppose that $P$ is an
$I$--bundle over a compact, orientable surface with non-empty boundary,
that the frontier $F$ of $P$ in $X$ is the vertical boundary of $P$,
and that $F$ is properly embedded in $X$.  Let $Y$ be a component of
$\overline{X{-}P}$, and assume that every component of
$\overline{X{-}(P\cup Y)}$ is a solid torus. (The last condition holds
vacuously if $\overline{X{-}P}$ is connected.) Then
$$\rank\pi_1(X)\le\rank\pi_1(Y)+\rank\pi_1(P).$$
\EndLemma

\Proof
Let $S$ denote the base of the $I$--bundle $P$.  Let $m+1$ denote the
number of components of $F$.  If $m=0$ then $F$ is connected, and the
assertion of the lemma is immediate from the Seifert-van Kampen
theorem. Hence we may assume that $m\ge1$.

Set $r=\rank\pi_1(Y)$ and $s=\rank\pi_1(P)=\rank\pi_1(S)$. We have
$s\ge m$, with equality if and only if $S$ is planar.

Let us denote the components of $\overline{X{-}P}$ by $Z_0,\ldots,Z_n$,
where $n\ge 0$, $Z_0=Y$, and $Z_j$ is a solid torus for each $j$ with
$0<j\le n$.  Let $A_0,\ldots,A_m$ denote the components of $F$, which
we index in such a way that $A_j\subset Z_j$ for $0\le j\le n$.  For
each $i\in\{0,\ldots,m\}$, the annulus $A_i$ is contained in a unique
component $Z_{q(i)}$ of $\overline{X{-}P}$; thus
$q:\{0,\ldots,m\}\to\{0,\ldots,n\}$ is a well-defined surjection, and
our indexing of the $A_i$ implies that $q(i)=i$ for $i=0,\ldots,n$.
(In particular $n\le m$.)

For $i=0,\ldots,m$, fix a point $a_i\in
A_i$, and fix a loop $\alpha_i$ in $A_i$, based at $a_i$, which
represents a generator of $\pi_1(A_i,a_i)$. Let $\beta_0$ denote the
constant path at $a_0$, and for each $i$ with $0<i\le m$ let $\beta_i$
denote a path in $P$ from $a_0$ to $a_i$, which projects to an
embedded arc $B_i$ in $S$. We may suppose the $\beta_i$ to be chosen
so that the arcs $B_i$ meet only at the point $a_0$. For
$i=0,\ldots,m$ let  $c_i$ denote the element of $\pi_1(P,a_0)$ represented by the
loop $\beta_i\ast\alpha_i\ast\overline{\beta_i}$. Then
$\pi_1(P,a_0)$ has a minimal generating set $\{x_0,\ldots,x_{s-1}\}$ such
that $x_i=c_i$ for $i=0,\ldots,m-1$.

For each $i$ with $n<i\le m$, we fix a path $\gamma_i$ in $Z_{q(i)}$
from $a_{q(i)}$ to $a_{i}$, and define an element $t_i\in\pi_1(X,a_0)$
by $t_i=[\beta_{q(i)}\ast\gamma_i\ast\overline{\beta_i}]$, where
brackets denote the based homotopy class of a loop in $X$. We fix a
generating set $\{y_1,\ldots,y_r\}$ for $\pi_1(Y,a_0)$. For
$j=1,\dots,n$ we fix  a loop $\zeta_j$  in $Z_j$ based at $a_j$ which
represents a generator for the cyclic group $\pi_1(Z_j,a_j)$, and we
set $z_j=[\beta_j\ast\zeta_j\ast\overline{\beta_j}]\in\pi_1(X,a_0)$. Then
$\pi_1(X,a_0)$ is generated by the set
$$\{\what{x}_i:0\le i<s\}\cup\{\what{y}_k:1\le k\le r\}\cup\{z_j:1\le j\le
n\}\cup\{t_i:n< i\le m\},$$
where $\what{x}_i$ and $\what{y}_k$ denote the images of 
${x_i}$ and $y_k$
under the inclusion homomorphisms from $\pi_1(P,a_0)$ and
$\pi_1(Y,a_0)$ to $\pi_1(X,a_0)$. 

For $i=0,\ldots,m$, let $ \what{c}_i$ denote the image of $c_i$
under the inclusion homomorphism $\pi_1(P,a_0)\to\pi_1(X,a_0)$.  Since
$\alpha_0$ is a loop in $A_0\subset Y$, and $\beta_0$ is the constant
path, $\what{c}_0$ is a word in the $\what{y}_k$.  Likewise, for
$0<i\le n$, since $\alpha_i$ is a loop in $A_i\subset Z_i$, the
element $\what{c}_i=[\beta_{i}\ast\alpha_i\ast\overline{\beta_i}]$
is a power of ${z_i}$. Furthermore, if $n<i\le m$, and if we set
$j=q(i)$, then $\what{c}_i=t_i^{-1}u_it_i$, where
$u_i=[\beta_{j}\ast\gamma_i\ast\alpha_i\ast\overline{\gamma_i}\ast\overline{\beta_j}]$.
Since $\gamma_i\ast\alpha_i\ast\overline{\gamma_i}$ is a loop in $Z_j$
based at $a_i$, it follows that $u_i$ is a word in the $\what{y}_k$
if $q(i)=0$, and is a power of $z_j$ if $q(i)>0$.

In particular,
the elements
$\what{c}_0,\ldots,\what{c}_m$ all lie in
the subgroup of $\pi_1(X,a_0)$ generated by the $\what{y}_i$, the $z_j$ and
the $t_i$. Since we have $\what{x}_i=\what{c}_i$ for
$i=0,\ldots,m-1$, we deduce that $\pi_1(X,a_0)$ is generated by the set
$$\{\what{x}_i:m\le i<s\}\cup\{\what{y}_k:1\le k\le r\}\cup\{z_j:1\le j\le
n\}\cup\{t_i:n< i\le m\}.$$
Hence 
$$\rank\pi_1(X)\le(s-m)+r+n+(m-n)=r+s,$$
which is the conclusion of
the lemma.
\EndProof  

The following lemma is a refined version of \cite[Lemma 2.21]{first}.

\Lemma\label{better 2.21} If $\calw$ is a connected normal book
of $I$--bundles, the rank of $\pi_1(|\calw|)$ is at most
 $2\barchi(|\calw|)+1$. 
\EndLemma

\Proof 
Set $W=|\calw|$. If $W$ is an $I$--bundle over a closed surface,
we have
$$\rank\pi_1(W)=\chibar(W)+2;$$
by normality we have $\chibar(W)>0$, so the conclusion holds in this
case. Hence we may assume that every page of $\calw$ meets at least
one binding.

Let $p$ denote the number of pages of $\calw$. We shall
recursively construct a finite sequence of connected sub-books
$\calw_1,\ldots,\calw_p$ of $\calw$, where $\calw_i$ has exactly $i$
pages. To begin the recursion, we choose an arbitrary page $P_1$ of
$\calw$ and define $\calw_1$ to consist of $P_1$ and the bindings that
meet it. Now assume that $\calw_i$ has been constructed for a given
$i<p$. Since $\calw $ has $p$ pages, $\calw_i$ is a proper sub-book of
$\calw$. Since $\calw $ is connected, $|\calw_i|$ must meet some page
$P_{i+1}$ not contained in $|\calw_i|$. We define $\calw_{i+1}$ to
consist of the pages and bindings of $\calw_i$, the page $P_{i+1}$, and
all bindings of $\calw$ that meet $P_{i+1}$.

We set $W_i=|\calw_i|$ for $i=1,\ldots,p$, and we let $W_0$ denote
some (arbitrarily chosen) binding of $\calw_1$. We shall show that
for $i=0,\ldots,p$ the rank of $\pi_1(W_i)$ is at most
 $2\barchi(W_i)+1$. For $i=0$ this is obvious, and for $i=p$ it is the
 conclusion of the lemma. It therefore suffices to show that if 
$0<k\le p$ and $\pi_1(W_{k-1})$ has rank at most $2\barchi(W_{k-1})+1$, then
$\pi_1(W_{k})$ has rank at most $2\barchi(W_{k})+1$.

The hypotheses of \fullref{new lemma} hold if we set $X=W_{k}$,
$P=P_{k}$, and $Y=W_{k-1}$. Hence we have
\begin{eqnarray*}
\rank\pi_1(W_{k})&\le\rank\pi_1(W_{k-1})+\rank\pi_1(P_{k})\\
&\le(2\chibar(W_{k-1})+1)+\rank\pi_1(P_{k}).
\end{eqnarray*}
Now since $P_{k}$ is an $I$--bundle over a compact, connected
surface with non-empty boundary, we have
$\rank\pi_1(P_{k})=1+\chibar(P_{k})$. But by the definition of a
normal book of $I$--bundles we have $\chibar(P_{k})\ge1$, and hence 
$$\rank\pi_1(P_{k})\le2\chibar(P_{k}).$$
It follows that
$$\rank\pi_1(W_{k})\le(2\barchi(W_{k-1})+1)+2\barchi(P_{k})=2\barchi(W_{k})+1,$$
as required.
\EndProof

Our next result is a analogue, in the context of the present section, of
\cite[Theorem 9.13]{first}.

\Lemma\label{better 9.14}
Assume that \fullref{my
    conjecture} is true.
Let $M$ be a closed simple $3$--manifold with $\rank\pi_1(M) \ge
11$. Suppose that $\dim H_1(M;\ZZ_2)\ge5$, and that
$\pi_1(M)$ has a
subgroup isomorphic to a genus--$2$ surface group. Then there
is a connected, normal book of $I$--bundles $\calw$ with
$W=|\calw|\subset M$ such that $\partial W$ is incompressible in
$M$ and $\barchi(W) = 2$.
\EndLemma

\Proof We shall adapt the proof of \cite[Theorem 9.13]{first}.  The
latter result has the same conclusion as the present lemma, but in
place of the hypothesis $\rank\pi_1(M)\ge11$, that $\dim
H_1(M;\ZZ_2)\ge5$, and that \fullref{my conjecture} is true, it
has the hypothesis that $\dim H_1(M;\ZZ_2)$ has rank at least $11$.
This hypothesis is used twice in the proof of \cite[Theorem
9.13]{first}: once in the first sentence to allow the application of
\cite[Corollary 9.11]{first}, and again in the fifth sentence of the
fifth and final paragraph of the proof. The application of
\cite[Corollary 9.11]{first} requires only the lower bound of $5$ for
$\dim H_1(M;\ZZ_2)$. Hence, under the hypotheses of the present lemma,
the first four paragraphs of the proof of \cite[Theorem 9.13]{first},
and the first sentence of the fifth paragraph, go through without
change, and show that either (a) there is a connected, normal book of
$I$--bundles $\calw$ with $W=|\calw|\subset M$ such that $\partial W$
is incompressible in $M$ and $\barchi(W) = 2$ (this is the case $m=0$
in the notation of the proof in \cite{first}) or (b) there exist a
finite-sheeted covering space $\tM$ of $M$ and a connected normal book
of $I$--bundles $\calw$ with $W=|\calw|\subset \tM$ such that the
inclusion homomorphism $\pi_1(W)\to \pi_1(\tM)$ is surjective and
$\chibar(W)\le 4$.  (Alternative (b) corresponds to the case $m>0$ in
the notation of the proof in \cite{first}, and we take $\tM=N_{m-1}$
in the notation of that proof. Since $N_{m-1}$ is closed, it is a
finite-sheeted covering space of $M$ according to \cite[Section~9.3]{first}.)

Now by \fullref{better 2.21} we have
$\rank\pi_1(|\calw|)\le2\barchi(|\calw|)+1\le9$.  The surjectivity of
$\pi_1(W)\to \pi_1(\tM)$ therefore implies that $\rank\pi_1(
\tM)\le9$.  On the other hand, since by hypothesis we have
$\rank\pi_1(M)\ge11$, \fullref{my conjecture}
implies that $\rank\pi_1( \tM)\ge 10$. This is a contradiction, and the
proof is complete.  \EndProof

The following result follows from \fullref{better 9.14} above in
exactly the same way that \cite[Corollary 9.14]{first} follows from
\cite[Theorem 9.13]{first}.

\Lemma\label{better 9.15} Assume that \fullref{my conjecture} is true. Let $M$ be a closed
simple $3$--manifold with $\rank\pi_1(M)\ge 11$. Suppose that
  $\pi_1(M)$ has a
  subgroup isomorphic to a genus--$2$ surface group, and that
  $\dim H_1(M;\ZZ_2)\ge5$.  Then $M$
contains either a connected incompressible surface of genus $2$ or a
separating, connected incompressible surface of genus $3$.\NoProof
\EndLemma

The following result is an analogue of \cite[Proposition 10.5]{first}
in the context of the present section.

\Lemma\label{better 9.5}
Assume that \fullref{my
  conjecture} is true. Suppose that $M$ is a closed orientable hyperbolic
3--manifold with $\rank\pi_1(M)\ge11$. Suppose that
$\pi_1(M)$ has a subgroup isomorphic to a genus--$2$ surface group, and
that $\dim H_1(M;\ZZ_2)\ge5$.
Then $\vol M\ge\WHAT$.
\EndLemma

\Proof
It follows from \fullref{better 9.15}  that either 
\begin{enumerate}
\item[(i)] $M$ contains either a separating incompressible surface of
genus $2$ or $3$, or
\item[(ii)] $M$ contains a non-separating incompressible surface of
genus $2$.
\end{enumerate}

Suppose that (i) holds but that { $\vol M<\WHAT$}.  Let $X_1$ and
$X_2$ denote the closures of the components of $M-S$. According to
\cite[Theorem 10.4]{first} (a result deduced from the main result of
Agol, Dunfield, Storm and Thurston \cite{ast}), each $X_i$ has the form $|\calw_i|$ for some book of
$I$--bundles $\calw_i$. For $i=1,2$ we have
$$\chibar(X_i)=\tfrac12\chibar(S)\le2.$$
By \fullref{better 2.21}, applied to $\calw=\calw_i$, it follows that
$$\rank\pi_1(X_i)\le2\barchi(X_i)+1\le5$$
for $i=1,2$.
Hence by the Seifert-van Kampen theorem, $$\rank \pi_1(M)\le\rank\pi_1( X_1)+\rank\pi_1( X_2)\le10,$$
a contradiction to the hypothesis. 

Now suppose that (ii) holds but that { $\vol M<\WHAT$}.  Let $X$
denote the connected manifold obtained by splitting $M$ along
$S$. According to \cite[Theorem 10.4]{first} we have $X=|\calw|$ for some
book of $I$--bundles $\calw$.  We have
$$\chibar(X)=\tfrac12\chibar(\partial X)=\chibar(S)=2.$$
By \fullref{better 2.21} it follows that
$$\rank\pi_1(X)\le2\barchi(X)+1\le5.$$
Hence $$\rank\pi_1 (M)\le\rank\pi_1 (X)+1\le6,$$
and again we have a contradiction. 
\EndProof

\Proof[Proof of \fullref{B}] Assume that the conclusion is
false, i.e. that $\rank\pi_1(M)\ge 11$ and that
$\dim H_1(M;\ZZ_2)\ge5$.  If $\pi_1(M)$ has a subgroup
isomorphic to a genus--$2$ surface group, then it follows from
\fullref{better 9.5} that $\vol M\ge\WHAT>3.08$, a
contradiction to the hypothesis.
  
  There remains the possibility that $\pi_1(M)$ has no subgroup
  isomorphic to a genus--$2$ surface group. Now since 
$H_1(M;\ZZ_2)$ has rank at least $5$, it follows from
\cite[Proposition 1.1]{shalen-wagreich}
 that every subgroup of rank at most $3$ in
  $\pi_1(M)$ has infinite index. But it follows from Anderson, Canary,
  Culler and Shalen \cite[Proposition 7.3,
  Remark 7.5]{accs} that if $M$ is an orientable hyperbolic $3$--manifold
  without cusps such that $\pi_1(M)$ contains no genus--$2$ surface
  subgroup and such that every subgroup of rank at most $3$ in
  $\pi_1(M)$ has infinite index, then $\pi_1(M)$ is $3$--free, in the
  sense that each subgroup of rank at most $3$ is free. And according
  to \cite[Corollary 10.3]{first}, if a closed, orientable hyperbolic
  $3$--manifold has $3$--free fundamental group, then its volume exceeds
  $3.08$. Again the hypothesis
  is contradicted.  \EndProof

Since Conjectures \ref{my conjecture} and \ref{genus version} are
equivalent modulo the ``rank equals genus'' \fullref{rank
    is genus},
\fullref{B} has the following immediate consequence.

\Corollary\label{BCor} If Conjectures \ref{rank is genus} and \ref{genus version} are
 true, then for every closed, orientable hyperbolic
$3$--manifold $M$ with volume at most $3.08$, either $M$ has Heegaard
 genus at most $10$, or
$\dim H_1(M;\ZZ_2)\le4$.  \EndCorollary

%
%
%
\bibliographystyle{gtart}
\bibliography{link}

\begin{thebibliography}{}
\providecommand\bibmarginpar{\leavevmode\marginpar}
\def\urlstyle#1{{\tt #1}}

\bibitem{agoltameness}
\textbf{I Agol}, \emph{Tameness of hyperbolic 3--manifolds}
  \xox{arXiv}{math.GT/0405568}

\bibitem{rank4}
\textbf{I Agol}, \textbf{M Culler}, \textbf{P Shalen}, \emph{Dehn surgery,
  homology and hyperbolic volume}  \xox{arXiv}{math.GT/0508208}

\bibitem{ast}
\textbf{I Agol}, \textbf{N\,M Dunfield}, \textbf{P\,A Storm}, \textbf{W\,P
  Thurston}, \emph{Lower bounds on volumes of hyperbolic Haken 3--manifolds}
  \xox{arXiv}{math.DG/0506338}

\bibitem{accs}
\textbf{J\,W Anderson}, \textbf{R\,D Canary}, \textbf{M Culler}, \textbf{P\,B
  Shalen}, \emph{Free {K}leinian groups and volumes of hyperbolic
  3--manifolds}, J. Differential Geom. 43 (1996) 738--782 \xox{MR}{1412683}

\bibitem{bz}
\textbf{M Boileau}, \textbf{H Zieschang}, \emph{Genre de Heegaard d'une
  vari\'et\'e de dimension 3 et g\'en\'erateurs de son groupe fondamental}, C.
  R. Acad. Sci. Paris S\'er. I Math. 296 (1983) 925--928 \xox{MR}{719278}

\bibitem{cg}
\textbf{D Calegari}, \textbf{D Gabai},
  \href{http://dx.doi.org/10.1090/S0894-0347-05-00513-8} {\emph{Shrinkwrapping
  and the taming of hyperbolic 3--manifolds}}, J. Amer. Math. Soc. 19 (2006)
  385--446 \xox{MR}{2188131}

\bibitem{second}
\textbf{M Culler}, \textbf{P\,B Shalen}, \emph{Singular surfaces, mod 2
  homology and hyperbolic volume II}  \xox{arXiv}{math.GT/0701666}

\bibitem{paradoxical}
\textbf{M Culler}, \textbf{P\,B Shalen},
  \href{http://dx.doi.org/10.2307/2152768} {\emph{Paradoxical decompositions,
  2--generator {K}leinian groups, and volumes of hyperbolic 3--manifolds}}, J.
  Amer. Math. Soc. 5 (1992) 231--288 \xox{MR}{1135928}

\bibitem{first}
\textbf{M Culler}, \textbf{P\,B Shalen}, \emph{Singular surfaces, mod 2
  homology and hyperbolic volume I}, Trans. Amer. Math. Soc.  (to appear)
  \xox{arXiv}{math.GT/0506396}

\bibitem{JS}
\textbf{W\,H Jaco}, \textbf{P\,B Shalen}, \emph{Seifert fibered spaces in
  3--manifolds}, Mem. Amer. Math. Soc. 21 (1979) viii+192 \xox{MR}{539411}

\bibitem{jones-reid}
\textbf{K\,N Jones}, \textbf{A\,W Reid},
  \href{http://dx.doi.org/10.1090/S0002-9939-00-05775-0} {\emph{Vol3 and other
  exceptional hyperbolic 3--manifolds}}, Proc. Amer. Math. Soc. 129 (2001)
  2175--2185 \xox{MR}{1825931}

\bibitem{mr}
\textbf{Y Moriah}, \textbf{H Rubinstein}, \emph{Heegaard structures of
  negatively curved 3--manifolds}, Comm. Anal. Geom. 5 (1997) 375--412
  \xox{MR}{1487722}

\bibitem{prezbetti}
\textbf{A Przeworski}, \href{http://dx.doi.org/10.1112/S0024609301008669}
  {\emph{Volumes of hyperbolic 3--manifolds of Betti number at least 3}}, Bull.
  London Math. Soc. 34 (2002) 359--360 \xox{MR}{1887708}

\bibitem{reid}
\textbf{A\,W Reid}, \emph{Some remarks on 2--generator hyperbolic
  3--manifolds}, from: ``Discrete groups and geometry (Birmingham, 1991)'',
  London Math. Soc. Lecture Note Ser. 173, Cambridge Univ. Press, Cambridge
  (1992)  209--219 \xox{MR}{1196929}

\bibitem{rr}
\textbf{Y Rieck}, \textbf{J\,H Rubinstein}, \emph{Invariant Heegaard surfaces
  in non-Haken manifolds with involutions} , preprint, University of Arkansas

\bibitem{rs}
\textbf{Y Rieck}, \textbf{E Sedgwick},
  \href{http://dx.doi.org/10.1016/S0166-8641(99)00147-9} {\emph{Persistence of
  {H}eegaard structures under {D}ehn filling}}, Topology Appl. 109 (2001)
  41--53 \xox{MR}{1804562}

\bibitem{schultweid}
\textbf{J Schultens}, \textbf{R Weidmann}, \emph{On the geometric and algebraic
  rank of graph manifolds}  (2003), preprint

\bibitem{shalen-wagreich}
\textbf{P\,B Shalen}, \textbf{P Wagreich},
  \href{http://dx.doi.org/10.2307/2154149} {\emph{Growth rates,
  $\mathbb{Z}_p$--homology, and volumes of hyperbolic 3--manifolds}}, Trans.
  Amer. Math. Soc. 331 (1992) 895--917 \xox{MR}{1156298}

\bibitem{census}
\textbf{J Weeks}, \textbf{C Hodgson}, \emph{Census of closed hyperbolic
  3--manifolds}
\ Available at \setbox0\hbox{\makeatletter\@url
{ftp://www.geometrygames.org/priv/weeks/SnapPea/SnapPeaCensus/ ClosedCensus/Clo%
sedCensusInvariants.txt}}
\href{ftp://www.geometrygames.org/priv/weeks/SnapPea/SnapPeaCensus/ClosedCensu%
s/ClosedCensusInvariants.txt}
{\unhbox0}

\end{thebibliography}

\end{document}